\documentclass[12pt]{amsart}
\usepackage{amscd}     
\usepackage{amssymb}
\usepackage{amsmath, amsthm, graphics, dsfont}
\usepackage{xypic}     
\LaTeXdiagrams        
\usepackage[all]{xy}
\xyoption{2cell} \UseAllTwocells \xyoption{frame} \CompileMatrices
\allowdisplaybreaks[3]
\usepackage{amsfonts}
\usepackage[normalem]{ulem}
\usepackage{hyperref}
\usepackage{tikz}
\usepackage{mathtools}
\usepackage{verbatim} 
\usepackage{MnSymbol}
\usepackage{latexsym}
\usepackage{epsfig}
\usepackage{enumerate}
\usepackage{times}
\usepackage{xcolor}
\usepackage{geometry}

\newgeometry{vmargin={25mm,25mm}, hmargin={25mm,25mm}}   % set the margins

\newtheorem{prop}{Proposition}[section]

\numberwithin{equation}{section}

\newcommand{\C}{\mathcal{C}}

\newcommand{\M}{\mathcal{M}}

\newcommand{\K}{\mathcal{K}}

\newcommand{\Z}{\mathcal{Z}}

\date{\today}

\begin{document}

\title{On counting plurifibered varieties}

\author{Yunfeng Jiang}
\address{Department of Mathematics\\ University of Kansas\\ 405 Snow Hall, 1460 Jayhawk Blvd. \\ Lawrence, KS 66045\\ USA}
\email{y.jiang@ku.edu}

\author{Hsian-Hua Tseng}
\address{Department of Mathematics\\ Ohio State University\\ 100 Math Tower, 231 West 18th Ave. \\ Columbus,  OH 43210\\ USA}
\email{hhtseng@math.ohio-state.edu}

\begin{abstract}
We consider the problem of enumeration of maps from plurifibered varieties. 
\end{abstract}

\maketitle

{\centering \em Dedicated to Dan Abramovich on the occasion of his 60th birthday\par}

\section{Introduction}

\subsection{Gromov-Witten theory}
Let $Y$ be a smooth projective variety over the complex number. Gromov-Witten theory of $Y$, which is an important aspect of modern enumerative geometry, is built from Kontsevich's moduli stacks of stable maps from nodal curves to $Y$. Gromov-Witten invariants of $Y$ may be viewed as virtual counts of maps from curves to $Y$. The term ``virtual'' here refers to the use of {\em virtual fundamental classes}. This has the disadvantage of making Gromov-Witten invariants usually not enumerative. On the other hand, Gromov-Witten invariants defined using virtual fundamental classes have very nice properties, e.g. the structure of a {\em cohomological field theory} \cite{km}.

The problem of extending Gromov-Witten theory to higher dimensional domains has been proposed for some time, see \cite{A2}. Compact (coarse) moduli spaces of stable maps from suitable higher dimensional domains have been constructed in \cite{A} and \cite{DR}. At this point, virtual fundamental classes for such moduli spaces of stable maps have not been constructed.
In the degree zero case, the moduli space is the product of the moduli space of varieties of general type and the target variety.  The KSBA moduli space of log general type varieties  is a very important topic in birational geometry, see \cite{Kollar-Shepherd-Barron}, \cite{HMX14}, \cite{KP15}. In the case of dimension two, the virtual fundamental class was constructed in \cite{Jiang_2022} for the moduli space of surfaces of general type. 

The purpose of this note is to discuss an extension of Gromov-Witten theory to a special kind of higher-dimensional domains, the so-called {\em plurifibered varieties} \cite{Ab}. Compact moduli stacks of maps from plurifibered varieties to $Y$ are described in Section \ref{sec:moduli}. In Section \ref{sec:obstruction}, we show that our moduli stacks can be equipped with obstruction theories. In general, these obstruction theories are {\em not} perfect. We discuss in Section \ref{sec:integral} two situations in which we do have perfect obstruction theories, hence virtual fundamental classes. Section \ref{sec:ci} contains a discussion on the case of complete intersections.

\subsection{Acknowledgment}

Both Y. Jiang and H.-H. T. are supported in part by Simons Foundation Collaboration Grant.

\section{Results}

\subsection{Moduli}\label{sec:moduli}
For a proper Deligne-Mumford stack $\mathcal{X}$, we use $\K_{g,n}(\mathcal{X},d)$ to denote the moduli stack of $n$-pointed genus $g$ degree $d$ twisted stable maps to $\mathcal{X}$ \cite{av}.

Fix a natural number $k$. Let $g_1,...,g_k\in \mathbb{Z}_{\geq 0}$, $n_1,...,n_k\in \mathbb{Z}_{\geq 0}$. For $d_1\in H_2(Y,\mathbb{Z})$, consider
\begin{equation}
\K_1(Y):=\K_{g_1,n_1}(Y, d_1).    
\end{equation}
We define inductively
\begin{equation}
\K_{i+1}(Y):= \K_{g_{i+1}, n_{i+1}}(\K_i, d_{i+1}), \quad \text{ for a choice of } d_{i+1}\in H_2(\K_i, \mathbb{Z}).    
\end{equation}
\iffalse
Put
\begin{equation}
\K_{\vec{g},\vec{n}}(Y, \vec{d}):=\K_{k}(Y).    
\end{equation}
\fi
By \cite{av}, $\K_i(Y)$, $i=1,...,k$ are proper Deligne-Mumford stacks with projective coarse moduli spaces. 

For a scheme $S$, an $S$-object of $\K_k(Y)$ is of the form
\begin{equation}
\xymatrix{
\C_1\ar[d]\ar[r] & \C_2\ar[d]\ar[r] &...\ar[r] & \C_{k-1}\ar[d]\ar[r] & \C_k\ar[d]\ar[r] & S\\ 
Y & \K_1(Y) & \, & \K_{k-2}(Y) & \K_{k-1}(Y). & \, 
}   
\end{equation}
Here for $i=1,...,k-1$, $\C_{i}\to \C_{i+1}$ is a family of twisted curves of genus $g_i$. The diagram 
\begin{equation*}
    \xymatrix{
\C_1\ar[r]\ar[d] & Y\\
\C_2, & \,
    }
\end{equation*}
is a $n_1$-pointed genus $g_1$ degree $d_1$ stable map to $Y$, and 
\begin{equation*}
\xymatrix{
\C_i\ar[r]\ar[d] & \K_{i-1}\\
\C_{i+1} & \,
}    
\end{equation*}
for $2\leq i\leq k$ are $n_i$-pointed genus $g_i$ degree $d_i$ stable maps to $\K_{i-1}$.

The sequence $$\bar{\C}_1\longrightarrow \bar{\C}_2\longrightarrow...\longrightarrow \bar{\C}_k\longrightarrow S$$ between coarse moduli spaces gives $\bar{\C}_1$ the structure of an $S$-family of {\em plurifibered varieties}, in the sense of \cite{Ab}. For this reason, we interpret $\K_k(Y)$ as a compact moduli stack of stable maps from plurifibered varieties to $Y$.

The marked gerbes on $\C_1,...,\C_k$ yield sections on coarse moduli spaces:
\begin{equation*}
\xymatrix{
\bar{\C}_1 \ar[r] & \bar{\C}_2\ar@/^1.5pc/[l]^{s^1_1,...,s^1_{n_1}} \ar[r] & \bar{\C}_3\ar@/^1.5pc/[l]^{s^2_1,...,s^2_{n_2}}\ar[r] & ......\ar@/^1.5pc/[l]^{s^3_1,...,s^3_{n_3}} \ar[r] & \bar{\C}_k\ar[r] \ar@/^1.5pc/[l]^{s^{k-1}_1,...,s^{k-1}_{n_{k-1}}} & S\ar@/^1.5pc/[l]^{s^k_1,...,s^k_{n_k}}.
}    
\end{equation*}
For $\vec{j}=(j_1, j_2,...,j_k)$ with $1\leq i_1\leq n_1$, $1\leq j_2\leq n_2$, ..., $1\leq j_k\leq n_k$, the composition 
$$s^1_{j_1}\circ...\circ s^{k-1}_{j_{k-1}}\circ s^k_{j_k}$$
with the map $\bar{\C}_1\to Y$ yields a map $S\to Y$. This construction yields an evaluation map
\begin{equation}\label{eqn:ev_map}
ev_{\vec{j}}: \K_k(Y) \longrightarrow Y.   
\end{equation}
We may use $ev_{\vec{j}}$ to pull back cohomology (or Chow) classes from $Y$ to $\K_k(Y)$.

\subsection{Obstruction theory}\label{sec:obstruction}
We consider the following general situation. Let $\mathcal{Z}$ be a proper Deligne-Mumford stack. Consider the universal family:
\begin{equation*}
\xymatrix{
\C\ar[d]^{\pi}\ar[r]^f & \mathcal{Z}\\
\K_{g,n}(\mathcal{Z},d). & \,
}    
\end{equation*}
Suppose $\phi: E^\bullet\to L_\Z$ is a perfect obstruction theory on $\Z$. Let $$T^{vir}_\mathcal{Z}$$ be the {\em virtual} tangent bundle associated to the perfect obstruction theory $\phi$ on $\mathcal{Z}$. By definition, $T^{vir}_\mathcal{Z}$ is the derived dual of $E^\bullet$. 
\begin{prop}\label{prop:obs_theory}
The object $R^\bullet\pi_*f^*T^{vir}_\mathcal{Z}$ is an obstruction theory.
\end{prop}

To simplify notations, we consider the Proposition in the case of fixed domain. The case of families of domains is similar.

Let $\C$ be an orbifold curve. Let $X=\text{Hom}(\C, \Z)$ be the Hom stack \cite{O}. Let $f: \C\times X\to \Z$, $\pi:\C\times X\to X$ be the universal family. We have maps 
\begin{equation*}
f^*L_\Z\to L_{\C\times X}\to L_{\C\times X/\C}\simeq \pi^*L_X,    
\end{equation*}
which gives a map $f^*L_\Z\to \pi^*L_X$. Composing with $f^*E^\bullet\to f^*L_\Z$, we get $f^*E^\bullet\to \pi^*L_X$. Tensoring with $\omega_\C$ yields
\begin{equation*}
    f^*E^\bullet\otimes^L \omega_\C\to \pi^*L_X\otimes^L\omega_\C=\pi^!L_X.
\end{equation*}
By duality, we get
\begin{equation*}
   E^{'\bullet}:=R\pi_*(f^*E^{\bullet\vee})^\vee \simeq R\pi_*(f^*E^\bullet\otimes^L\omega_\C)\to L_X.
\end{equation*}
We claim that this is an obstruction theory. To see this, we use \cite[Theorem 4.5]{BF}. 

Let $T\to \bar{T}$ be a square zero extension with ideal sheaf $J$. Consider $g: T\to X$. An extension $$\bar{g}: \bar{T}\to X$$ exists if and only if the obstruction class, which is a class in $$Ext_T^1(g^*L_X, J),$$ vanishes. Since $X=\text{Hom}(\C, \Z)$, a map $g: T\to X$ is equivalent to $h: \C\times T\to \Z$. An extension $\bar{g}: \bar{T}\to X$ of $g$ is equivalent to $$\bar{h}: \C\times \bar{T}\to \Z$$ which restricts to $h$. Since $\C\times T\to \C\times \bar{T}$ is a square zero extension with ideal sheaf $p^*J$, where we write $p: \C\times T\to T$ for the projection, the existence of $\bar{h}$ is equivalent to the vanishing of the obstruction class, which is a class in $$Ext_{\C\times T}^1(h^*L_\Z, p^*J).$$ Since $E^\bullet\to L_\Z$ is an obstruction theory, the existence of $\bar{h}$ is equivalent to the vanishing of the image of the obstruction class in $Ext_{\C\times T}^1(h^*E^\bullet, p^*J)$ (\cite[Theorem 4.5]{BF}).

Consider the cartesian diagram:
\begin{displaymath}
    \xymatrix{ 
    \mathcal{C}\times T\ar[r]_{id\times g}\ar[d]_{p}\ar@/^1.5pc/[rr]^{h} & \C\times X\ar[d]^{\pi}\ar[r]_{f} & \Z \\
   T\ar[r]_{g} & X. &}
\end{displaymath}
Apply \cite[Lemma 6.1]{BF} to this diagram, we have for any $k$
\begin{equation*}
Ext_{\C\times T}^k(h^*E^\bullet, p^*J)=Ext_{\C\times T}^k((id\times g)^*f^*E^\bullet, p^*J)\simeq Ext_T^k(g^*\underbrace{R\pi_*(f^*E^\bullet\otimes \omega_\C)}_{E^{'\bullet}}, J).    
\end{equation*}
Therefore, the existence of $\bar{g}$, which is equivalent to the existence of $\bar{h}$, is equivalent to the vanishing of the image of the obstruction class in $$Ext_T^1(g^*E^{'\bullet}, J).$$ 

A similar argument shows that extensions of $g$ form a torsor under $Ext_T^0(g^*E^{'\bullet}, J)$. Therefore by \cite[Theorem 4.5]{BF} $E^{'\bullet}\to L_X$ is an obstruction theory on $X$.

\subsection{Integrations}\label{sec:integral}
Applying Proposition \ref{prop:obs_theory} repeatedly, we obtain an obstruction theory on $\K_k(Y)$. %$\K_{\vec{g},\vec{n}}(Y,\vec{d})$. 
In general, this obstruction theory is {\em not} perfect. Here we discuss two cases in which we do get perfect obstruction theories.
\subsubsection{Surfaces}\label{sec:surface}
Consider the $k=2$ case:
\begin{equation}\label{eqn:moduli_fibered_maps}
\K_{2}(Y)=\K_{g_2,n_2}(\K_{g_1,n_1}(Y,d_1),d_2).
\end{equation}
Objects of this moduli stack are maps to $Y$ from {\em fibered} surfaces, c.f. \cite{AbV}. 

If $g_1=0$ and $Y$ is convex (i.e. $H^1(C, f^*TY)=0$ for any genus $0$ stable map $f:C\to Y$), then $\K_{g_1,n_1}(Y,d_1)$ is a smooth DM stack and hence (\ref{eqn:moduli_fibered_maps}) has a perfect obstruction theory by \cite{agv}. Therefore (\ref{eqn:moduli_fibered_maps}) has a virtual fundamental class $$[\K_2(Y)]^{vir},$$ 
which is a class in homology (or Chow) group of $\K_2(Y)$.  

We can be a bit more explicit about the virtual dimension of $\K_2(Y)$, as follows. The evaluation maps 
\begin{equation*}
ev_j: \K_{2}(Y)=\K_{g_2,n_2}(\K_{0,n_1}(Y,d_1),d_2)\to \bar{I}\K_{0,n_1}(Y, d_1), \quad j=1,...,n_2,    
\end{equation*}
take value in the rigidified inertia stack $\bar{I}\K_{0,n_1}(Y,d_1)$. For connected components
\begin{equation*}
\mathcal{I}_1,...,\mathcal{I}_{n_2}\subset \bar{I}\K_{0,n_1}(Y,d_1),    
\end{equation*}
the virtual dimension of 
$$\bigcap_{j=1}^{n_2}ev_j^{-1}(\mathcal{I}_j)\subset \K_{g_2,n_2}(\K_{0,n_1}(Y,d_1),d_2)$$ is given by 
\begin{equation*}
(1-g_2)\text{ dim }\K_{0,n_1}(Y,d_1)+\int_{d_2}c_1(T_{\K_{0,n_1}(Y,d_1)}) +3g_2-3+n_2-\sum_{j=1}^{n_2}\text{ age }(\mathcal{I}_j).     
\end{equation*}
Here 
\begin{equation*}
\text{ dim }\K_{0,n_1}(Y,d_1)=\text{ dim }Y +\int_{d_1}c_1(T_Y)+n_1-3,    
\end{equation*}
and a formula for the tangent bundle $T_{\K_{0,n_1}(Y,d_1)}$ can be found in e.g. \cite[Section 2.5.3]{c_thesis}. This formula is based on the virtual dimension formula for moduli stacks of twisted stable maps to Deligne-Mumford stacks, see \cite{agv}.

We can consider integrals of the form:
\begin{equation}\label{eqn:integral}
 \int_{[\K_2(Y)]^{vir}}\prod_{\vec{j}}ev_{\vec{j}}^* \alpha_{\vec{j}}\in \mathbb{Q}, \quad \alpha_{\vec{j}}\in H^*(Y). \end{equation}
The classes $ev_{\vec{j}}^* \alpha_{\vec{j}}$ can be viewed as imposing incidence conditions given by Poincar\'e duals of $\alpha_{\vec{j}}$. (\ref{eqn:integral}) can be interpreted as virtual counts of ruled surfaces in $Y$. 

Examples of convex varieties are projective spaces, Grassmannians, and other homogeneous spaces.

In the special case $Y=\text{point}$, $\K_{g_1,n_1}(\text{point},0)$ is smooth for any $g_1$. Hence $\K_2(\text{point})$ always admits a perfect obstruction theory. Integrations against the virtual fundamental class of $\K_2(\text{point})$ are Gromov-Witten invariants of moduli stacks of stable pointed curves, which is an interesting subject on its own (see e.g. \cite{ms}).

\subsection{Complete intersections}\label{sec:ci}

We discuss a heuristic situation. Let $Z$ be a smooth projective variety, $\mathcal{V}\to Z$ a vector bundle over $Z$, and $s: Z\to \mathcal{V}$ a regular section. Consider the zero locus 
$$Y:=s^{-1}(0).$$ 
Let $\M(-)$ be a moduli space of stable maps from some kind of domains. Then $$\M(Y)\subset \M(Z)$$ is the zero locus of the section $\tilde{s}$ of the {\em sheaf} $R^0\pi_*f^*\mathcal{V}$ induced from $s$. Here 
\begin{equation*}
\xymatrix{
\mathcal{U}\ar[d]^{\pi}\ar[r]^{f} & Z\\
\mathcal{M} & \,
}    
\end{equation*}
is the universal family. In view of twisted Gromov-Witten theory \cite{cg}, we can define 
\begin{equation}\label{eqn:vir_class_formula}
    [\M(Y)]^{tw}:=[\M(Z)]^{vir}\cap e_{\mathbb{C}^*}(R^\bullet\pi_*f^*\mathcal{V}).
\end{equation}
Here $e_{\mathbb{C}^*}(-)$ is the $\mathbb{C}^*$-equivariant Euler class. The object $R^\bullet \pi_*f^*\mathcal{V}$ admits a $\mathbb{C}^*$ action induced from the $\mathbb{C}^*$ action on $\mathcal{V}$ by scaling fibers.

The existence of the object $R^\bullet\pi_*f^*\mathcal{V}$ in $K$-theory requires some properties of the universal family $\pi:\mathcal{U}\to \M$, namely $\M$ needs to has resolution property and $\pi$ needs to be lci.

(\ref{eqn:vir_class_formula}) allows us to define $[\M(Y)]^{tw}$ if $[\M(Z)]^{vir}$ exists. As considered in Section \ref{sec:surface}, if $Z$ is convex and $\M$ is $\K_{g_2,n_2}(\K_{0,n_1}(-,d_1),d_2)$, then $[\M(Z)]^{vir}$ exists. By \cite{grpan}, $\K_{0,n_1}(Z,d_1)$ is a quotient stack. Then by \cite{agot}, $\M(Z)=\K_{g_2,n_2}(\K_{0,n_1}(Z,d_1),d_2)$ is also a quotient stack and has resolution property. The universal family $\pi$ can be described as follows:
\begin{equation*}
\xymatrix{
\C_1\ar[r]_{\pi_2'}\ar[d]\ar@/^1.5pc/[rr]^{\pi} & \C_2\ar[r]_(.25){\pi_1}\ar[d]^{f_1} & \M(Z)=\K_{g_2,n_2}(\K_{0,n_1}(Z,d_1),d_2)\\
Z & \K_{0,n_1}(Z,d_1). & \,
}    
\end{equation*}
Here $\pi_1:\C_2\to \K_{g_2,n_2}(\K_{0,n_1}(Z,d_1),d_2)$ is the universal family and $\pi_2'$ is the pullback by $f_1$ of the universal family over $\K_{0,n_1}(Z,d_1)$. 

By \cite{FaP}, $\pi_2'$ is lci. By \cite{agot}, $\pi_1$ is lci. Hence $\pi$ is lci. Therefore (\ref{eqn:vir_class_formula}) applies.


\begin{thebibliography}{12}

\bibitem{Ab} D. Abramovich, {\em Canonical models and stable reduction for plurifibered varieties}, arXiv:math/0207004. 

%\bibitem{AbC} D. Abramovich, J.-C. Chen, {\em Flops, flips and perverse point sheaves on threefold stacks}, J. Algebra 290 (2005), 372--407.

\bibitem{agot} D. Abramovich, T. Graber, M. Olsson, H.-H. Tseng, {\em On The Global Quotient Structure of The Space of Twisted Stable Maps to a Quotient Stack}, J. Alg. Geom. 16 (2007),731--751.

\bibitem{agv} D. Abramovich, T. Graber, A. Vistoli, {\em Gromov-Witten theory of Deligne-Mumford stacks}, Amer. J. of Math. 130 (2008), no. 5, 1337--1398.


%\bibitem{AbH} D. Abramovich, B. Hassett, {\em Stable varieties with a twist}, in: ``Classification of algebraic varieties'', 1--38, EMS Ser. Congr. Rep., Eur. Math. Soc., Z\"urich, 2011.

\bibitem{AbV} D. Abramovich, A. Vistoli, {\em Complete moduli for fibered surfaces}, in: ``Recent progress in intersection theory (Bologna, 1997)'', 1--31, Trends Math., Birkh\"auser Boston, Boston, MA, 2000. 

\bibitem{av} D. Abramovich, A. Vistoli, {\em Compactifying the space of stable maps}, J. Amer. Math. Soc. 15 (2002), 27--75.


%\bibitem{ABIP} K. Ascher, D. Bejleri, G. Inchiostro, Z. Patakfalvi, {\em Wall crossing for moduli of stable log pairs}, arXiv:2108.07402, to appear in Ann. of Math. 

\bibitem{A} V. Alexeev, {\em Moduli spaces $M_{g,n}(W)$ for surfaces}, in: ``Higher-dimensional complex varieties (Trento, 1994)'', 1--22, de Gruyter, Berlin (1996).

\bibitem{A2} V. Alexeev, {\em Higher-dimensional analogues of stable curves}, in: ``Proceedings of the international congress of mathematicians (ICM), Madrid, Spain, August 22–30, 2006. Volume II: Invited lectures'', 515--536, European Mathematical Society.

%\bibitem{AG} V. Alexeev, G. M. Guy, {\em Moduli of weighted stable maps and their gravitational descendants}, J. Inst. Math. Jussieu 7 (2008), no. 3, 425--456. 

\bibitem{BF} K. Behrend, B. Fantechi, {\em The intrinsic normal cone}, Invent. Math. 128 (1997), no. 1, 45--88. 

%\bibitem{BeI} D. Bejleri, G. Inchiostro, {\em Stable pairs with a twist and gluing morphisms for moduli of surfaces}, Selecta Math. (N.S.) 27 (2021), no. 3, Paper No. 40, 44 pp.

%\bibitem{B} C. Birkar, {\em Moduli of algebraic varieties}, arXiv:2211.11237.

\bibitem{c_thesis} T. Coates, {\em Riemann-Roch theorems in Gromov-Witten theory}, Ph.D. thesis, UC Berkeley, 2003. 

\bibitem{cg} T. Coates, A. Givental, {\em Quantum Riemann-Roch, Lefschetz and Serre}, Ann. of Math. (2) 165 (2007), no. 1, 15--53.


\bibitem{DR} R. Dervan, J. Ross, {\em Stable maps in higher dimensions}, Math. Ann. 374 (2019), 1033--1073.

\bibitem{FaP} C. Faber, R. Pandharipande, {\em Hodge integrals and Gromov-Witten theory}, Invent. Math. 139, No. 1, 173--199 (2000).

%\bibitem{FP} W. Fulton, R. Pandharipande, 

\bibitem{grpan} T. Graber, R. Pandharipande, {\em Localization of virtual classes}, Invent. Math. 135 (1999), 487--518.

\bibitem{HMX14} C. D. Hacon, J. MacKernan, and C. Xu, \newblock Boundedness of moduli of varieties of general type, J. Euro. Math. Soc. 20 (2018), Issue 4, 865-901,   arXiv:1412.1186. 

\bibitem{Jiang_2022}Y. Jiang, \newblock The virtual fundamental class for the moduli space of surfaces of general type, arXiv:2206.00575.

%\bibitem{K} J. Koll\'ar, {\em Families of varieties of general type}, Cambridge Tracts in Mathematics, 231. Cambridge University Press, Cambridge, 2023. xviii+471 pp, also available on the author's webpage.

\bibitem{Kollar-Shepherd-Barron} J. Kollar and N. I. Shepherd-Barron, \newblock Threefolds and deformations of surface singularities, 
{\em Invent. Math.}, 91, 299-338 (1998).

\bibitem{km} M. Kontsevich, Y. Manin, {\em Gromov-Witten classes, quantum cohomology, and enumerative geometry}, Comm. Math. Phys. 164 (1994), 525--562.


\bibitem{KP15} S. Kov\'acs, Z. Patakfalvi, {\em Projectivity of the moduli space of stable log-varieties and subadditivity of log-Kodaira dimension}, J. Amer. Math. Soc. 30 (2017), no. 4, 959--1021. 

\bibitem{ms} Y. Manin, M. Smirnov, {\em Towards motivic quantum cohomology of $\bar{M}_{0,S}$}, Proc. Edinb. Math. Soc., II. Ser. 57, No. 1, 201--230 (2014), arXiv:1107.4915.

\bibitem{O} M. Olsson, {\em Hom-stacks and restriction of scalars}, Duke Math. J. 134 (2006), no. 1, 139--164.

%\bibitem{O2} M. Olsson, {\em A boundedness theorem for Hom-stacks}, Math. Res. Lett. 14 (2007), no. 6, 1009--1021.



\end{thebibliography}
\end{document}